\documentclass[letterpaper,10pt,conference,twoside]{ieeeconf}

\IEEEoverridecommandlockouts                              

\usepackage{amsmath,amssymb,color}
\usepackage{stackengine}
\usepackage{graphicx}

\usepackage{latexsym}
\usepackage{multicol}
\usepackage{algorithm}
\usepackage{algpseudocode}

\usepackage{verbatim}
\usepackage{cite}
\usepackage[english]{babel}
\usepackage[latin1]{inputenc}
\usepackage{enumerate}

\usepackage{amsmath,amsthm}
\usepackage{bbm}

\usepackage{pgf}
\usepackage{pgfplots,tikz}
\usetikzlibrary{tikzmark}
\usepackage[font=footnotesize]{caption}
\usepackage{arydshln}

\usepackage{amsfonts}
\usepackage{amsbsy}
\usepackage{bm}

\newcommand\oprocendsymbol{\hbox{$\bullet$}}
\newcommand\oprocend{\relax\ifmmode\else\unskip\hfill\fi\oprocendsymbol}

\DeclareMathOperator{\conv}{conv}

\DeclareMathOperator*{\argmin}{arg\,min}

\DeclareMathOperator{\interior}{int}
\DeclareMathOperator{\uniform}{Uni}

\let\leq\leqslant
\let\geq\geqslant

\newcommand{\R}{\mathbb R}
\newcommand{\N}{\mathbb N}

\newcommand{\calF}{\ensuremath{\mathcal{F}}}

\newcommand{\calS}{\ensuremath{\mathcal{S}}}

\newcommand{\calZ}{\ensuremath{\mathcal{Z}}}
\newcommand{\lse}{\operatorname{lse}}

\newcounter{todocounter}
\usepackage[colorinlistoftodos]{todonotes}

\newtheorem{theorem}{Theorem}[section]

\newtheorem{lemma}[theorem]{Lemma}
\newtheorem{corollary}[theorem]{Corollary}

\newtheorem{remark}[theorem]{Remark}

\newtheorem{assumption}{Assumption}

\newtheorem{problem}{Problem}

\allowdisplaybreaks
\renewcommand{\baselinestretch}{0.985}

\title{Set-valued regression and cautious suboptimization: From noisy
  data to optimality\thanks{This work was partially funded by ONR
    Award N00014-23-1-2353 and SNF/FW Weave Project 200021E\_20397}}

\author{Jaap Eising \quad Jorge Cort\'{e}s\thanks{Jaap Eising is with
    the Department of Information Technology and Electrical
    Engineering at ETH Z\"{u}rich, Switzerland,
    \texttt{jeising@ethz.ch}.  Jorge Cort\'{e}s is with the Department
    of Mechanical and Aerospace Engineering, University of California,
    San Diego, \texttt{cortes@ucsd.edu}. The bulk of this work was
    done when J. Eising was affiliated with UC San Diego.}}

\begin{document}

\maketitle
	
\begin{abstract}
  This paper deals with the problem of finding suboptimal values of an
  unknown function on the basis of measured data corrupted by bounded
  noise. As a prior, we assume that the unknown function is
  parameterized in terms of a number of basis functions. Inspired by
  the informativity approach, we view the problem as the
  suboptimization of the worst-case estimate of the function. The
  paper provides closed form solutions and convexity results for this
  function, which enables us to solve the problem. After this, an online
  implementation is investigated, where we iteratively measure the
  function and perform a suboptimization. This nets a procedure that
  is safe at each step, and which, under mild assumptions, converges
  to the true optimizer.
\end{abstract}

\section{Introduction}

Optimization of unknown functions on the basis of measured data is a
topic with many applications ranging from control to machine
learning. For instance, cost or reward functions in modern
model-predictive control methods are often partially unknown due to
modeling difficulties. On the other hand, sampling the functions is
often possible, giving us access to measurements. Accurately
determining suboptimal values of such unknown functions is at best a
major part of some control objective and at worst
safety-critical. This motivates the problem: determine suboptimal
points of the unknown function on the basis of noisy measurements. We
will investigate this problem both in a \emph{one-shot} setting, where
the data are given, and in an \emph{online} setting that allows for
repeated measurements.

\emph{Literature review:} Of course this is not a new problem, and
there are many solutions with various setups of which we mention the
ones most relevant. In order to solve the one-shot problem, the usual 
approach is to first employ the data to obtain a unique or, in some sense,
`best' estimate of the unknown function. For this, a number of
nonparametric techniques have been developed to estimate unknown
functions from data. Popular are e.g. Gaussian processes
\cite{CER-CKIW:05} and methods based on Lipschitz constants
\cite{MM-AV:91, MM-CN:02,JC-SJR-CER-JM:20}. However, in this paper, we
will make the assumption that the unknown function can be
parameterized in terms of a number of basis functions or
features. Common choices of basis functions are for instance linear,
polynomial, Gaussian, or sigmoidal functions. These basis
functions allow us to perform regression (see
e.g. \cite{TH-RT-JF:13,RT:96}) on the parameters, leading to the
`best' estimate. Methods differ on what is considered this best
estimate: For instance least squares (minimal Frobenius norm), ridge
regression (minimal $L_2$ norm), sparse or Lasso regression (minimal
$L_1$ norm). Analogously
are the similar methods that have arisen in the nonlinear
system-identification literature. Methods have been developed to 
determine models that are sparse \cite{SLB-JLP-JNK:16}, low rank
\cite{PJS:10,PJS:22,JNK-SLB-BWB-JLP:16}, or both \cite{MRJ-PJS-JWN:14} 
within a class parameterized using a basis.

After obtaining an estimate of the unknown function, one can treat
this estimate as the true unknown function and apply any well-studied
optimization technique to obtain suboptimal values. Given the fact
that the data is corrupted by noise, it is reasonable to require some
robustness from the methods applied (see 
e.g. \cite{AB-LEG-AN:09,DB-DBB-CC:11} and the references
therein). Apart from these one-shot optimization problems, we are also
interested in online problems, in which we iteratively measure and
optimize. This is inspired by methods such as extremum-seeking control
\cite{MK-HHW:00,KBA-MK:03,ART-DP:01}, whereas our implementation is
markedly different.

In contrast to this paradigm of regression (or: learning) followed by 
optimization,  this paper can be viewed as being in line with the concept
of informativity (see \cite{HJVW-JE-HLT-MKC:20,HJVW-MKC-JE-HLT:22}),
where system properties are investigated for \emph{all} systems compatible
with the measurements. This falls within a recent surge of
replacing system identification with methods on the basis of Willems'
fundamental lemma \cite{JCW-PR-IM-BLMDM:05}. While most of these works
deal with linear dynamics, the paper \cite{IM:22} dealt with bilinear
systems by embedding them into a higher dimensional linear system. 
Similarly, the works \cite{MG-CDP-PT:21,CDP-MR-PT:23} consider
systems that are linear in their basis functions. 

\emph{Statement of contributions:} As mentioned, we take a viewpoint
related to that of the informativity framework. Using the assumed
basis functions, we characterize the set of all parameters compatible
with a set of measurements with bounded noise. This is named
\emph{set-valued regression}. Clearly, a point can be guaranteed to be
suboptimal for the unknown function only if it is for \emph{all}
functions corresponding to compatible parameters. This motivates us to
introduce \emph{cautious suboptimization}, that is, the problem of finding
such values in a way that is robust against the worst-case realization of 
the parameters. In order to resolve this, we derive closed forms and
investigate convexity of this realization. Combining these
allows us to resolve the problem using any method from convex
optimization. Importantly, this allows us to derive \emph{guaranteed}
upper bounds of the optimal value of the unknown functions on the
basis of potentially very small data sets.

In addition to this one-shot problem, we investigate an online variant
consisting of the iteration of two steps: Collecting local
measurements and performing cautious suboptimization. This gives rise
to a procedure which given increasingly sharp upper bounds of the
unknown function. Moreover, assuming that the noise is randomly
generated in addition to being bounded, we prove that this procedure
converges to the true optimal value.

Proofs are omitted for reasons of space and will appear elsewhere.
\section{Problem statement}
Let\footnote{Throughout the paper, we use the following notation. We
  denote by $\N$ and $\R$ the sets of nonnegative integer and real
  numbers, respectively. We let $\R^{n\times m}$ denote the space of
  $n\times m$ real matrices. For vectors $v\in\R^n$, we write $v\geq0$
  (resp. $v>0$) for elementwise nonnegativity (resp. positivity). The
  sets of such vectors are denoted
  $\mathbb{R}^n_{\geq0}:=\{v\in\R^n | v\geq0\}$ and
  $\mathbb{R}^n_{>0}:=\{v\in\R^n | v>0 \}$. On the other hand, for
  $P\in\R^{n\times n}$, $P\geq 0$ (resp. $P> 0$) denotes that $P$ is 
  symmetric positive semi-definite (resp. definite). We denote the 
  smallest singular value of $M\in\mathbb{R}^{n\times m}$ by 
  $\sigma_-(M)$. For a set 
  $\calS\subseteq\mathbb{R}^n$ we denote the convex hull by 
  $\conv(\calS)$ and the interior by $\interior(\calS)$.}
$\phi_i:\mathbb{R}^n\rightarrow\mathbb{R}$ for $i=1,\ldots, k$ be a
collection of known \emph{basis functions} (or \emph{features}) and
consider the set consisting of all functions
$\phi^\gamma:\mathbb{R}^n\rightarrow\mathbb{R}$ linearly parameterized by $\gamma\in\mathbb{R}^k$ as
$ \phi^\gamma(z) =
  \Sigma_{i=1}^{k}\gamma_i\phi_i(z)$.
By collecting the basis functions in a
vector-function as
\[
  b(z) := \begin{bmatrix} \phi_1(z) & \cdots &
    \phi_k(z)\end{bmatrix}^\top,
\]
we can write the shorthand $\phi^\gamma(z) = \gamma^\top b(z)$.
Consider a function $\hat{\phi}:\mathbb{R}^n\rightarrow\mathbb{R}$
which is unknown but can be expressed as a linear combination of the
features, i.e., $\hat{\phi}(z)=\phi^{\hat{\gamma}}(z)$, for some
\textit{unknown} parameter $\hat{\gamma}$.

We are interested in deducing properties of the function $\hat{\phi}$
on the basis of measurements. Suppose we sample the function for the
variables $z_i$ with $i=1,\ldots, T$ and collect noisy measurements
$y_i$ of the true function, that is, $ y_i =\hat{\phi}(z_i) +
w_i$. Here, the vector $w_i$ denotes an unknown noise, or disturbance,
for each $i$. In order to reason with these measurements in a
structured manner, we define 
\begin{align}\label{eq:Y-W}
  Y :=
  \begin{bmatrix}
    y_1 & \cdots
    & y_T
  \end{bmatrix},
      \quad W
      :=
      \begin{bmatrix}
        w_1 & \cdots & w_T
      \end{bmatrix},  
\end{align}
\[
  \Phi := \begin{bmatrix} b(z_1) & \cdots & b(z_T)\end{bmatrix}
  = \scalebox{0.8}{$\begin{bmatrix} \phi_1(z_1) & \dots & \phi_1(z_T) \\ \vdots & &
    \vdots \\ \phi_k(z_1) & \dots & \phi_k(z_T) \end{bmatrix}$}.
\]
Note that $Y,W\in \mathbb{R}^{1\times T}$,
$\Phi \in \mathbb{R}^{T\times k}$, and $Y = \hat{\gamma}^\top\Phi +W$.

We consider bounded noise. 
In order to formalize this, let $\Pi\in\mathbb{R}^{(1+\ell)\times (1+\ell)}$
be a symmetric matrix. We partition $\Pi$ as
\[
\Pi =
\begin{bmatrix}
	\Pi_{11} & \Pi_{12} \\ \Pi_{21} &
	\Pi_{22}
\end{bmatrix}, \textrm{ with } \Pi_{11}\in\mathbb{R}, \Pi_{22}\in \mathbb{R}^{\ell\times \ell}.
\] 
If $\Pi_{22}<0$, then we denote
the Schur complement
$\Pi|\Pi_{22}:= \Pi_{11}-\Pi_{12}\Pi_{22}^{-1}\Pi_{21}$. We define the set
\[
  \calZ(\Pi) := \left\lbrace v\in \mathbb{R}^\ell \mid \scalebox{0.8}{$\begin{bmatrix} 1
      \\ v\end{bmatrix}^\top \Pi \begin{bmatrix} 1 \\
      v \end{bmatrix} $} \geq 0 \right\rbrace .
\]
 We make the
following assumption on the noise model.

\begin{assumption}[Noise model]\label{as:noise model}
  Let $\Pi\in\mathbb{R}^{(1+T)\times(1+T)}$ be symmetric, such that
  $\Pi_{22}<0$, and $\Pi|\Pi_{22}\geq 0$.  The noise samples satisfy
  $W^\top\in \calZ(\Pi)$.
\end{assumption}

With $\Pi$ as in Assumption~\ref{as:noise model}, the set $\calZ(\Pi)$
is nonempty, convex, and bounded. A common example of such a noise
model is the case where $W W^\top \leq q$, for some $q\geq 0$, or as confidence intervals of Gaussian noise.
Assuming that the noise signal satisfies Assumption~\ref{as:noise
  model}, we can define the set of all parameters $\gamma$ consistent
with the measurements by:
\begin{equation}\label{eq:def Gamma}
  \Gamma := \{ \gamma \in \R^k \mid Y =
  \gamma^\top\Phi +W, W^\top\in\calZ(\Pi) \}.
\end{equation}
Thus, if we define $N\in\mathbb{R}^{(1+k)\times (1+k)}$ by
\begin{equation}\label{eq:defN}
  N:= \scalebox{0.9}{$\begin{bmatrix} N_{11} & N_{12} \\ N_{21} & N_{22}\end{bmatrix} $}
  = \scalebox{0.9}{$\begin{bmatrix} 1 & Y \\ 0& -\Phi \end{bmatrix}
  \Pi \begin{bmatrix} 1 & Y \\ 0& -\Phi \end{bmatrix}^\top$},  
\end{equation}
it follows immediately that $\Gamma = \calZ(N)$. Note that
$\hat{\gamma}\in\Gamma$ and that we have no further information on the
value of $\hat{\gamma}$. We refer to the procedure of obtaining
$\Gamma$ from the measurements as \emph{set-valued regression}.

\begin{remark}[Sufficiently exciting measurements]\label{rem:exciting}
  {\rm Note that the set $\Gamma$ is closed. Moreover, it is bounded if and only if
    $N_{22}<0$. Since $N_{22} = \Phi\Pi_{22}\Phi^\top$ and
    $\Pi_{22}<0$, this holds if and only if $\Phi$ has full row
    rank. In turn, this requires that the basis functions are not
    identical and that the set of points $z_i$ is `rich' enough or
    sufficiently `exciting', cf. \cite{JCW-PR-IM-BLMDM:05}. \oprocend
  }
\end{remark}

\begin{remark}[Least-squares estimates]\label{rem:lse}
  {\rm One can check that, if $N_{22}<0$, then
    $\gamma^{\lse}:=-N_{22}^{-1}N_{21} \in \Gamma$. Therefore $ (\gamma^{\lse})^\top b(z)$ is consistent with the
    measurements. In fact,
    \[
      \scalebox{0.9}{$\begin{bmatrix} 1 \\ -N_{22}^{-1}N_{21} \end{bmatrix}^\top
      N \begin{bmatrix} 1 \\
        -N_{22}^{-1}N_{21} \end{bmatrix} $} \geq \scalebox{0.9}{$\begin{bmatrix} 1 \\
        \gamma\end{bmatrix}^\top N \begin{bmatrix} 1 \\
        \gamma \end{bmatrix} $},
    \]
    for any $\gamma^\top\in \calZ(N)$. As such, $\gamma^{\lse}$ is the
    value for which the quadratic inequality is maximal. This leads us
    to refer to the function
    $\phi^{\lse}(z;\Gamma):=(\gamma^{\lse})^\top b(z) =
    -N_{12}N_{22}^{-1}b(z)$ as the \textit{least-squares estimate} of
    $\hat{\phi}(z)$. \oprocend }
\end{remark}

We are interested in the optimization of the unknown
function~$\hat{\phi}$.  However, based on the measurements, we cannot
distinguish between the different functions $\phi^\gamma$ for
$\gamma \in \Gamma$.  Indeed, small changes in the parameter $\gamma$
might lead to large changes in the quantitative behavior and the
location of optimal values of the functions $\phi^\gamma$.  In order
to be robust against such changes, we consider suboptimization
problems instead: for instance, we can conclude that
$\hat{\phi}(z)\leq \delta$ only if $\phi^\gamma(z) \leq \delta$ for all
$\gamma\in\Gamma$. This motivates the following.

\begin{problem}[Cautious optimization for set-valued
  regression]\label{prob:cautious-opt}
{\rm  Consider an unknown function $\hat{\phi}$, a noise model $\Pi$ such that Assumption~\ref{as:noise
    model} holds, measurements of the true function $(Y,\Phi)$, and
  $\Gamma$ as in \eqref{eq:def Gamma}. Then,
  \begin{enumerate}[(a)]
  \item\label{prob:verif}  (Verification of
      suboptimality): given $z\in\mathbb{R}^n$, find the smallest of $\delta\in\mathbb{R}$ for which
      $\phi^\gamma(z) \leq \delta$ for all $\gamma\in\Gamma$;
    \item\label{prob:cautious} (One-shot cautious suboptimization):
      using the solution to \eqref{prob:verif}, and given a set
      $\calS\subseteq \mathbb{R}^n$, find $z\in\calS$ for which
      \eqref{prob:verif} yields the minimal value of $\delta$;
    \item\label{prob:online} (Online cautious suboptimization):
      determine where to collect new measurements to iteratively
      improve the bound obtained in~\eqref{prob:cautious}.
    \end{enumerate} }
\end{problem}

Note that Problem~\ref{prob:cautious-opt} can
be posed instead as a question regarding properties of the measurements
$(Y,\Phi)$, as in the data informativity framework,
e.g.~\cite{HJVW-JE-HLT-MKC:20,HJVW-MKC-JE-HLT:22}. For instance, given
$\delta\in\mathbb{R}$ and $\calS\subseteq\mathbb{R}^n$, one could say
that the data $(Y,\Phi)$ is \emph{informative for
  $\delta$-suboptimization on $\calS$} if there exists $z\in\calS$
such that $\phi^{\gamma}(z)\leq \delta$ for all $\gamma\in\Gamma$.

\section{One-shot cautious suboptimization}
This section addresses Problems
\ref{prob:cautious-opt}.\eqref{prob:verif} and
\ref{prob:cautious-opt}.\eqref{prob:cautious}.  Consider 
\begin{equation}\label{eq:def g h}
  \phi^+(z;\Gamma) := \sup_{\gamma\in\Gamma} \phi^\gamma(z),\quad  \hspace{1em} \phi^-(z;\Gamma) := \inf_{\gamma\in\Gamma}
  \phi^\gamma(z),
\end{equation} 
which correspond to the elementwise worst-case realization of the
unknown parameter~$\hat{\gamma}$. Note that if $\Gamma$ is compact,
cf. Remark~\ref{rem:exciting}, the supremum and infimum are both
attained over $\Gamma$. Hence they can be replaced by maximum and
minimum, respectively, which implies that both $\phi^+(z;\Gamma)$ and
$\phi^-(z;\Gamma)$ are finite-valued functions. 

Resolving Problem~\ref{prob:cautious-opt}.\eqref{prob:verif} is
equivalent to determining function values of
$\phi^+(\cdot;\Gamma)$. Similarly, we can reformulate
Problem~\ref{prob:cautious-opt}.\eqref{prob:cautious} as finding
\begin{equation}\label{eq:minmax g} 
  \min_{z\in\calS} \max_{\gamma\in\Gamma} \phi^\gamma(z) =
  \min_{z\in\calS} \phi^+(z;\Gamma).  
\end{equation}
This problem takes the form of a minimax or bilevel optimization
problem. In this section we first investigate the inner problem of
finding values of $\phi^+(z;\Gamma)$, i.e., solving
Problem~\ref{prob:cautious-opt}.\eqref{prob:verif}.  After a detour
regarding \textit{uncertainty}, we check this function for
convexity. Then, by explicitly finding gradients, we can efficiently
resolve Problem~\ref{prob:cautious-opt}.\eqref{prob:cautious}.

\subsection{Verification of suboptimality}\label{sec:verification}

As a first step towards resolving the cautious suboptimization problem
we investigate \textit{verification} of suboptimality. That is,
\emph{given} $z\in\mathbb{R}^n$, test whether the unknown function is
such that $\hat{\phi}(z) \leq \delta$.  Recall that this problem can
be resolved if we can explicitly find function values of the functions
in~\eqref{eq:def g h}.  The following result provides closed-form
expressions for these functions on the basis of measurements.

\begin{theorem}[Closed-form expressions for bounds]
	\label{thm:explicit delta}
  Assume that the measurements $(Y,\Phi)$ are such that
  $\Gamma=\calZ(N)$ with $N_{22}<0$. Then
  \begin{align*}
    \phi^\pm(z;\Gamma)
    &=\! -N_{12}N_{22}^{-1}b(z) \pm
      \sqrt{(N|N_{22})b(z)^\top (-N_{22}^{-1}) b(z)}.
  \end{align*}
\end{theorem} 

Note that the first term in the closed forms expressions of the
functions $\phi^+(\cdot;\Gamma)$ and $\phi^-(\cdot;\Gamma)$ is the
least squares estimate, cf. Remark~\ref{rem:lse}. The result in
Theorem~\ref{thm:explicit delta} then shows that the difference
between the values of true unknown function and the least squares
estimate can be quantified in terms of the basis functions and the
data, as expressed in~$N$.

As a consequence of Theorem~\ref{thm:explicit delta}, we have the
following result expressing the gradient of the bounding functions. 

\begin{corollary}[Gradients in terms of data]\label{cor:gradient}
  Assume that the measurements $(Y,\Phi)$ are such that
  $\Gamma=\calZ(N)$ with $N_{22}<0$. Let the basis functions
  $\phi_i$ be differentiable and such that $b(z)\neq0$ for all
  $z\in\calS$.  Then,
  \begin{align*}
    \nabla \phi^\pm(z;\Gamma)
    & = [
      \nabla \phi_1(z) \cdots \nabla \phi_k(z) ] \cdot
    \\
    &   \scalebox{1}{$\bigg(
      \!-\! N_{22}^{-1}N_{21} \pm \sqrt{N|N_{22}}\frac{(-N_{22}^{-1})
      b(z)}{\sqrt{b(z)^\top (-N_{22}^{-1}) b(z)}}\bigg)$}.
  \end{align*}
\end{corollary}	
\subsection{Uncertainty of function values}\label{sec:uncertainty}
The discussion in Section~\ref{sec:verification} allows us to find
bounds for the unknown function $\hat{\phi}$, but does not consider
how much these bounds deviate from its true value.  By definition, we
have
$
  \phi^-(z;\Gamma) \leq \hat{\phi}(z) \leq \phi^+(z;\Gamma).
$
Therefore, we define the \textit{uncertainty at} $z$ by
\[
  U(z;\Gamma) := \phi^+(z;\Gamma) - \phi^-(z;\Gamma)
\]
to quantify how well we know the function value of $\hat{\phi}(z)$ at
$z\in\mathbb{R}^n$.  If the uncertainty at $z$ is close to 0, then the
function value of $\hat{\phi}(z)$ is quantifiably close to
$\phi^{\lse}(z;\Gamma) = -N_{12}N_{22}^{-1}b(z)$. To balance the
demands of a low upper bound on the value of the true function with an
associated low uncertainty, it is reasonable to consider the following
generalization of the cautious suboptimization
problem~\eqref{eq:minmax g}: for $\lambda\geq 0$, consider
\begin{equation}\label{eq:weighted cs}
  \min_{z\in\calS} \phi^+(z;\Gamma) +\lambda U(z;\Gamma).
\end{equation}
Under the conditions of Theorem~\ref{thm:explicit delta}, we obtain
the following closed form for the objective function.

\begin{lemma}[Explicit forms and~\eqref{eq:weighted
    cs}]\label{lem:expl uncertainty}
  Let the measurements $(Y,\Phi)$ be such that
  $\Gamma=\calZ(N)$ with $N_{22}<0$. Then
 $   U(z;\Gamma) = 2\sqrt{(N|N_{22}) b(z)^\top(- N_{22}^{-1}) b(z)}.$
  Moreover, if
  \[
    N_\lambda := \scalebox{0.9}{$\begin{bmatrix} N_{11} & N_{12} \\ N_{21} &
      N_{22} \end{bmatrix} $} + \scalebox{0.9}{$\begin{bmatrix}
      4\lambda(1+\lambda)(N|N_{22}) & 0\\ 0&0 \end{bmatrix}$},
  \]
  and $\Gamma_\lambda := \calZ(N_\lambda)$, then
$
    \phi^+(z;\Gamma) +\lambda U(z;\Gamma) =
    \phi^+(z;\Gamma_\lambda).
$
\end{lemma}

Lemma~\ref{lem:expl uncertainty} means that, even though
problem~\eqref{eq:weighted cs} is more general than the cautious
suboptimization problem~\eqref{eq:minmax g}, both problems can be
resolved in the same fashion.

\subsection{Convexity and suboptimization}
To provide efficient solutions to \eqref{eq:minmax g}, we investigate
when $\phi^+(\cdot;\Gamma)$ is convex. Towards this, we first
investigate conditions under which we can guarantee that the true
function $\hat{\phi}$ is convex.  Since nonnegative combinations of
convex functions are convex, the following result identifies
conditions that ensure the set of parameters consistent with the
measurements are nonnegative.

\begin{lemma}[Test for nonnegativity of parameters]\label{lem:gam
    positive}
  Given measurements $(Y,\Phi)$, let $\Gamma=\calZ(N)$, where $N$
  is as in \eqref{eq:defN}. Then,
  $\Gamma \subseteq \mathbb{R}^k_{\geq 0}$ if and only if $\Phi$ has
  full row rank
  and one of the following conditions hold
  \begin{enumerate}
  \item\label{lemitem:single case} The matrix $N\leq 0$ and
    $-N_{22}^{-1} N_{21}\in\mathbb{R}^k_{\geq 0}$, or
  \item\label{lemitem:scalar} The matrix $N$ has one positive
    eigenvalue, $-N_{22}^{-1} N_{21} \in \mathbb{R}^k_{>0}$ and, for
    all $i=1,\ldots,k$,
    \[
      N|N_{22} + \scalebox{1}{$\frac{(e_i^\top N_{22}^{-1} N_{21})^2}{e_i^\top
        N_{22}^{-1} e_i}$}\leq 0.
    \]
  \end{enumerate}
\end{lemma} 

We can use this result to identify conditions that ensure the
convexity of the unknown function and its upper bound.

\begin{corollary}[Convexity of the true function and upper
  bound]\label{cor:convexity}
  Suppose that the basis functions $\phi_i$ are convex and
  $\Gamma\subseteq \mathbb{R}_{\geq 0}^k$.
  \begin{itemize}
  \item Then, $\phi^\gamma$ is convex for all $\gamma\in\Gamma$ and
    $\phi^+(\cdot;\Gamma)$ is a finite-valued convex function;
  \item If, in addition, the functions $\phi_i$ are strictly convex
    and $0\not\in\Gamma$, then $\phi^\gamma$ is strictly convex for
    all $\gamma\in\Gamma$ and $\phi^+(\cdot;\Gamma)$ is strictly
    convex.
  \end{itemize}
\end{corollary}

Lemma~\ref{lem:gam positive} and Corollary~\ref{cor:convexity}
taken together mean that, if the basis functions are convex, we can
test for convexity on the basis of data.

Recall that we are interested in the optimization
problem \eqref{eq:minmax g}, and therefore not necessarily in
properties of the true function $\hat{\phi}$, but of its upper bound
$\phi^+(\cdot;\Gamma)$. This motivates our ensuing discussion to
provide conditions that ensure convexity of the upper bound instead.
Note that, under the assumptions of Theorem~\ref{thm:explicit delta},
we have
\[
  \phi^+(z;\Gamma) =\phi^{\lse}(z;\Gamma) + \tfrac{1}{2}
  U(z;\Gamma).
\]
Thus, if (i) $\phi^{\lse}(\cdot;\Gamma)= -N_{12}N_{22}^{-1}b(\cdot)$
is convex and (ii) $U(\cdot;\Gamma)$ is convex, then so is
$\phi^+(\cdot;\Gamma)$. Moreover, if in addition either is strictly convex, then
so is $\phi^+(\cdot;\Gamma)$.  Condition (i) could be checked directly
if all basis functions $\phi_i$ are twice continuously differentiable
by computing the Hessian of $\phi^{\lse}$.
Here, we present the following simple criterion derived from
composition rules, see e.g. \cite[Example 3.14]{SB-LV:09}, to test
for condition (ii).

\begin{corollary}[Convexity of the uncertainty]\label{cor:conv uncert}
  Assume that the measurements $(Y,\Phi)$ are such that
  $\Gamma=\calZ(N)$ with $N_{22}<0$. Then $U(\cdot;\Gamma)$ is
  convex if each basis function $\phi_i$ is convex and
  $-N_{22}^{-1}b(z)\geq 0$ for all $z\in\mathbb{R}^n$.
\end{corollary} 

%
%

Equipped with the results of this section, one can solve the cautious
suboptimization problems~\eqref{eq:minmax g} and~\eqref{eq:weighted
  cs} efficiently.  Under the assumptions of Theorem~\ref{thm:explicit
  delta}, we can write
closed-form expressions for $\phi^+(\cdot;\Gamma)$. This, in turn,
allows us to test for (strict) convexity using e.g.,
Corollaries~\ref{cor:convexity} or~\ref{cor:conv uncert}. If so, we
can apply (projected) gradient descent, using
Corollary~\ref{cor:gradient}, to resolve cautious suboptimization.

%


\section{Online cautious optimization}\label{sec:online}

In this section, we develop an online optimization procedure on the
basis of local measurements of the true function to refine the
optimality gap. Specifically, we devise a procedure where we first
collect data \textit{near} a candidate optimizer, we update a convex
upper bound of $\hat{\phi}$ on the basis of the measurements, and
lastly we recompute the candidate optimizer on the basis of the
updated upper bound.

To formalize this, we require some notation.  Let
$\calF =\{f_i\}_{i=1}^T \subseteq\mathbb{R}^n$ be a finite set. 
For a given $z\in\mathbb{R}^n$ we
measure the function at all points in $z +\calF$. For this, define
\[
  \Phi^\calF(z) := 
\scalebox{0.8}{$\begin{bmatrix} \phi_1(z+f_1) & \dots & \phi_1(z+f_T) \\ \vdots & &
  \vdots \\ \phi_k(z+f_1) & \dots & \phi_k(z+f_T) \end{bmatrix}$}.
\]
Given an initial point $z_0$, consider measurements at step~$k$,
\begin{equation}\label{eq:measure}
  Y_k = \hat{\gamma}^\top\Phi^\calF(z_k) + W_k, \textrm{ with }
  W_k^\top \in \calZ(\Pi),
\end{equation} 
where $Y_k$ and $W_k$ are as in~\eqref{eq:Y-W}. Define the set
$\Gamma_k$ of parameters which are compatible with the
$k^{\textup{th}}$ set of measurements,
\begin{equation}\label{eq:def Gk}
  \Gamma_k:=
  \calZ(N_k),
\end{equation}
where $
  N_k:= \scalebox{0.9}{$\begin{bmatrix} 1 & Y_k \\ 0& -\Phi^\calF(z_k) \end{bmatrix}
  \Pi \begin{bmatrix} 1 & Y_k \\ 0&
    -\Phi^\calF(z_k) \end{bmatrix}^\top$}$.

The online optimization procedure then incrementally incorporates
these measurements to refine the computation of the candidate
optimizer.  The following result investigates the properties of the
resulting \textit{online gradient descent}.

\begin{theorem}[Online gradient descent]\label{thm:online}
  Let $\calF$ be a finite set such that $0\in\interior (\conv \calF)$
  and such that $\Phi^\calF(z)$ has full row rank for all $z$. Define
  $\calS(z) := z +\conv\calF$. Consider an initial point
  $z_0\in\mathbb{R}^n$ and suppose that the measurements
  $(Y_0,\Phi^\calF(z_0))$ are such that $\phi^\gamma$ is strictly
  convex for all $\gamma\in\Gamma_0$. For $k\geq 1$, repeat
  the following two steps iteratively: first update the candidate
  optimizer
  \begin{equation}\label{eq:update est}
    z_k := \argmin_{z\in
      \calS(z_{k-1})} \phi^+(z;\Gamma_0\cap\ldots\cap\Gamma_{k-1})
  \end{equation} 
  and secondly measure the function $\hat{\phi}$ as in
  \eqref{eq:measure} and define $\Gamma_{k}$ as in~\eqref{eq:def Gk}.
  Then the following hold:
  \begin{enumerate} 
  \item\label{item:convex} For any $k\geq 1$, the problem
    \eqref{eq:update est} is strictly convex;
  \item\label{item:ubounds} For each $k\geq 1$, the algorithm provides
    an upper bound
    \begin{equation}\label{eq:ubounds}
      \min_{z\in\mathbb{R}^n}
      \hat{\phi}(z)\leq 
      \phi^+(z_k;\Gamma_0\cap\ldots\cap\Gamma_{k-1});
    \end{equation}
  \item\label{item:monotone} The upper bounds are monotonically
    nonincreasing
    \begin{equation}\label{eq:monotone}
      \phi^+(z_{k+1};\Gamma_0\cap\ldots\cap\Gamma_{k})  \leq
      \phi^+(z_k;\Gamma_0\cap\ldots\cap\Gamma_{k-1})  
    \end{equation}
  \item\label{item:strict mono} If $z_k\neq z_{k+1}$, then
    \eqref{eq:monotone} holds with a strict inequality.
  \end{enumerate} 
\end{theorem}

The algorithm described in Theorem~\ref{thm:online} provides a
sequence of upper bounds to true function values,
cf.~\eqref{eq:ubounds}, on the basis of \textit{local}
measurements. This means that after \textit{any} number of iterations,
we obtain a `worst-case' estimate of the function value
$\hat{\phi}(z_k)$ and, as such, for the minimum of $\hat{\phi}$. Given
that this sequence of upper bounds is nonincreasing,
cf.~\eqref{eq:monotone}, and bounded below by the true minimum of
$\hat{\phi}$, we can conclude that the algorithm converges. However,
without further assumptions, one cannot guarantee convergence of the
upper bounds to the minimal value of $\hat{\phi}$, or respectively of
$z_k$ to the global minimum of $\hat{\phi}$ (the simulations of
Section~\ref{sec:sims} below show an example of this precisely).

The following result shows that if the uncertainty is sufficiently
small near the optimizer, then the optimizer of the upper bound is
close to the global optimizer of the true, unknown function.

\begin{lemma}[Stopping criterion]\label{lem:stop}
  Let $\Gamma$ be compact and $\calS\subseteq\mathbb{R}^n$
  closed. 
  Define
  \[
    \bar{z}:=\argmin_{z\in\calS} \phi^+(z;\Gamma), \quad  
    \hat{z}:=\argmin_{z\in\calS} \hat{\phi}(z).
  \]
  Then
  $ \phi^+(\bar{z};\Gamma) \geq \hat{\phi}(\hat{z})\geq
  \phi^+(\bar{z};\Gamma) - \max_{z\in\calS} U(z;\Gamma)$.
\end{lemma} 

From Lemma~\ref{lem:stop}, we see that if the uncertainty on $\calS$
is equal to zero, then the local minima of $\phi^+$ and $\hat{\phi}$
coincide. In addition, if $\hat{\phi}$ is strictly convex, we have
that any local minimum in the \emph{interior} of $\calS$ is equal to
its global minimum.

When repeatedly collecting measurements, it seems reasonable to assume
that the uncertainty would decrease. However, without making further
assumptions, this is not necessarily the case.  In particular, a
situation might arise where repeated measurements corresponding to a
worst-case, noise signal give rise to convergence to a
fixed bound with nonzero uncertainty. To address this problem, we
consider a scenario where the noise samples are not only bounded but
distributed uniformly over the set $\calZ(\Pi)$.  To make this formal,
suppose that $\Pi$ is such that $\calZ(\Pi)$ is bounded. Consider the
measure $\mu$ and probability distribution over $\mathbb{R}^{T}$,
\[
  p(W) :=
  \scalebox{0.9}{$\begin{cases}
    \frac{1}{\mu(\calZ(\Pi))} &
    W^\top\in\calZ(\Pi), \\ 0 & \textrm{otherwise}.
  \end{cases}$}
\]
As a notational shorthand, we write
$W^\top\sim\uniform(\calZ(\Pi))$. The following result shows that,
under uniformly distributed noise samples, the uncertainty does indeed
decrease.

\begin{theorem}[Uncertainty under repeated
  measurements]\label{thm:stochastics}
  Under the assumptions of Theorem~\ref{thm:online}, suppose in
  addition that for $k\geq 1$, the measurements in \eqref{eq:measure}
  are such that $W_k^\top\sim \uniform(\calZ(\Pi))$ and
  $\sigma_-(\Phi^{\calF}(z))\geq a$ for all $z$. Then for any
  $z\in\mathbb{R}^n$, the expected value of the uncertainty
  monotonically converges to 0, that is,
  \[
    U(z;\Gamma_0\cap\ldots\cap \Gamma_{k-1})\geq
    U(z;\Gamma_0\cap\ldots\cap \Gamma_k),
  \]
  and
$
    \lim_{k\rightarrow \infty} \mathbb{E}(U(z;\Gamma_0\cap\ldots\cap
    \Gamma_k)) = 0.
$
\end{theorem}

As a consequence of Lemma~\ref{lem:stop} and
Theorem~\ref{thm:stochastics}, one can conclude that the expected
difference between the optimal value $\min_{z\in\calS} \hat{\phi}(z)$
of the unknown function and the optimal value
$\min_{z\in\calS} \phi^+(z;\Gamma_0 \cap \dots \cap \Gamma_k)$
provided by online gradient descent both converge to zero.


\section{Simulation examples}\label{sec:sims}
We illustrate here our results in a simple example. Let the unknown
function $\hat{\phi}:\mathbb{R}^2 \rightarrow \mathbb{R}$ be given by
$\hat{\phi}(z)=1+z^\top z$. Let $z=\begin{pmatrix} z_1& z_2\end{pmatrix}^{\!\top}$ and consider the basis functions
\[\phi_1(z) = 1, \hspace{1em}\phi_2(z) = z_1,  \hspace{1em}\phi_3(z) = z_2,  \hspace{1em}\phi_4(z) = z^\top z.\]

The value of the true parameter $\hat{\gamma}$ is
$\begin{pmatrix} 1 & 0& 0 & 1\end{pmatrix}^\top$.  We sample the
function at points in $z + \calF$, where
$\calF = \lbrace (0,0),(1,0),(0,1),(-1,-1)\rbrace$ (i.e., we measure
at the point itself and three points around it). For this choice,
$\Phi^\calF(z)$ has full row rank for all $z\in\mathbb{R}^2$.  The
measurements are corrupted by a large amount of noise: we assume a
noise model of the form $W_k W_k^\top \leq 30 $ for all $k\geq 0$ and
that the noise is uniformly distributed in this set. This means that
\[
  \Pi = \scalebox{0.9}{$\begin{bmatrix} 30 & 0\\ 0 & -I_4\end{bmatrix}$} \textrm{ and }
  W_k^\top \sim \uniform (\calZ(\Pi)).
\]
For $z_0$, we collect uniform random noisy measurements
$(Y_0,\Phi^\calF(z_0))$, leading to a set of consistent parameters
$\Gamma_0$. We take $z_0 = \begin{pmatrix} 3 & 3\end{pmatrix}^\top$
and verify that $\Gamma_0$ is such that $\phi^\gamma$ is strictly
convex for all $\gamma\in\Gamma_0$. 
%
On the basis of this, define
\[
  z_1 := \argmin_{z\in \calS(z_0)} \phi^+(z;\Gamma_0).
\]
Moreover, we can determine the least-squares estimate of the parameter
$\hat{\gamma}$ and the function $\hat{\phi}$. The latter is given by
\[
 \scalebox{0.9}{$ \phi^{\lse}(z;\Gamma_0)\! =\! 63.99\phi_1(z) - 23.27 \phi_2(z)
  -23.27\phi_3(z) + 5.1\phi_4(z)$}.
\]
Now we can evaluate the true function, least-squares estimate, and the
upper bound, finding the following values:
\[
  \scalebox{0.85}{$\begin{array}{llll}
  	\hline
    &z_0\quad\quad	&z_1\quad\quad	& 0\quad\quad
    \\
    \hline
    \hat{\phi}(\cdot)			&19 	&13.48 	& 1
    \\
    \phi^{\lse}(\cdot;\Gamma_0)\quad	&16.26 	&11.44	& 63.99
    \\
    \phi^+(\cdot;\Gamma_0)		&21.74	&16.74 	& 136.45
    \\
    \hline
  \end{array}$}
\]
We can make a few observations. First, the least-squares
estimate is quite close at the measured point $z_0$, yet far at the
true minimum of $\hat{\phi}$ at the origin. Further, the upper bound
indeed decreases monotonically. Lastly, while the upper
bound majorizes $\hat{\phi}$, the least-squares estimate does not.

\begin{figure}[ht]
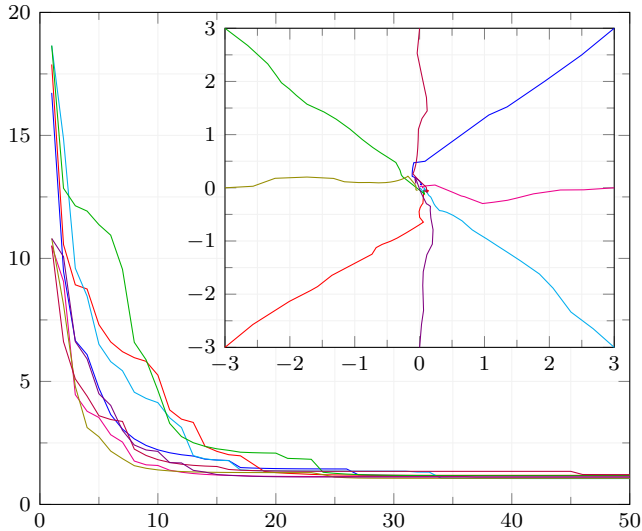

  \begin{center}
    \stackinset{r}{8pt}{t}{6pt}{\input{./tex/randomplottraj.tex}}{\input{./tex/randomplotbounds.tex}}
  \end{center}
 \vspace{-1em} \caption{The simulation results for eight different initial
    conditions. In the large plot we see the upper bounds
    corresponding to the different trajectories, which indeed decrease
    monotonically towards the true minimum $\hat{\phi}(0)=1$. The
    inset image shows each of the corresponding trajectories of the
    estimate $z_k$, revealing that each tends to the true minimizer.}
\label{fig:plots more meas}
\end{figure}

We run 100 iterations of the online gradient descent algorithm in
Theorem~\ref{thm:online} and plot the results in Figure~\ref{fig:plots
  single meas}. In spite of the relatively noisy data, we see rapid
convergence of the parameters~$\gamma$ and the values of the
estimate~$z_k$. Moreover, the upper bound is close to the true value
at the measured point. After just 10 steps, we see that the true
function value at the estimate is already significantly lower than at
the start. To show that this happens regardless of the choice of
initial conditions, Figure~\ref{fig:plots more meas} shows the
trajectories of $z_k$ resulting from 8 different initial conditions
and the corresponding values of the upper bounds.

\pgfplotsset{
	compat=1.11,
	legend image code/.code={
		\draw[mark repeat=2,mark phase=2,thick]
		plot coordinates {
			(0cm,0cm)
			(0.15cm,0cm)        
			(0.3cm,0cm)         
		};%
	}
}
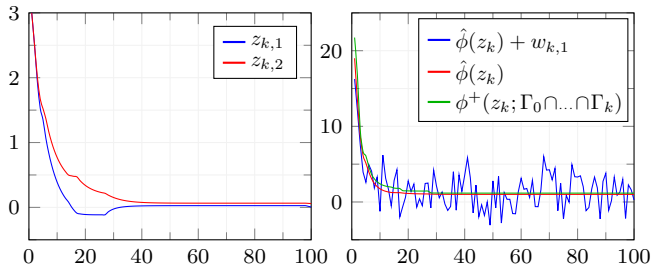
\begin{figure}[ht]
\begin{center}
\vspace{-1em}\hspace{-0.4em}\begin{tikzpicture}
	\begin{axis}[%
		width=0.3\textwidth,
		xmin=0,
		xmax=100,
		ymin=-0.5,
		ymax=3,
		grid=both,
		grid style={line width=.1pt, draw=gray!10},
		minor tick num=1,
		tick label style={font=\footnotesize},
		legend cell align=left,
		legend style={at={(axis cs:95,2.8)},anchor=north east,nodes={scale=0.8, transform shape}}
		]
		\addplot[smooth, color=blue] table[row sep=crcr]{%
			1	3       \\
			2	2.498203\\
			3	1.931903\\
			4	1.531562\\
			5	1.349320\\
			6	1.067195\\
			7	0.843843\\
			8	0.660895\\
			9	0.511041\\
			10	0.388294\\
			11	0.287752\\
			12	0.205397\\
			13	0.137940\\
			14	0.084185\\
			15	0.012644\\
			16	-0.04565\\
			17	-0.08961\\
			18	-0.09961\\
			19	-0.10588\\
			20	-0.10975\\
			21	-0.11205\\
			22	-0.11338\\
			23	-0.11408\\
			24	-0.11441\\
			25	-0.11451\\
			26	-0.11448\\
			27	-0.11438\\
			28	-0.07281\\
			29	-0.04461\\
			30	-0.02515\\
			31	-0.01149\\
			32	-0.00177\\
			33	0.005243\\
			34	0.010370\\
			35	0.014162\\
			36	0.016998\\
			37	0.019141\\
			38	0.020778\\
			39	0.022039\\
			40	0.023019\\
			41	0.023784\\
			42	0.024387\\
			43	0.024864\\
			44	0.025242\\
			45	0.025545\\
			46	0.025787\\
			47	0.025982\\
			48	0.026139\\
			49	0.026265\\
			50	0.026367\\
			51	0.026450\\
			52	0.026516\\
			53	0.026570\\
			54	0.026613\\
			55	0.026649\\
			56	0.026679\\
			57	0.026701\\
			58	0.026719\\
			59	0.026735\\
			60	0.026747\\
			61	0.026758\\
			62	0.026766\\
			63	0.026773\\
			64	0.026779\\
			65	0.026784\\
			66	0.026788\\
			67	0.026791\\
			68	0.026793\\
			69	0.026795\\
			70	0.026797\\
			71	0.026799\\
			72	0.026800\\
			73	0.026801\\
			74	0.026801\\
			75	0.026801\\
			76	0.026801\\
			77	0.026801\\
			78	0.026801\\
			79	0.026801\\
			80	0.026802\\
			81	0.026803\\
			82	0.026802\\
			83	0.026802\\
			84	0.026802\\
			85	0.026803\\
			86	0.026803\\
			87	0.026803\\
			88	0.026802\\
			89	0.026802\\
			90	0.026802\\
			91	0.026802\\
			92	0.026803\\
			93	0.026803\\
			94	0.026803\\
			95	0.026804\\
			96	0.026804\\
			97	0.026803\\
			98	0.026803\\
			99	0.016170\\
			100	0.007476\\
		};
		\addlegendentry{$z_{k,1}$}
		\addplot[smooth, color=red] table[row sep=crcr]{%
			1	3       \\
			2	2.498203\\
			3	2.006030\\
			4	1.672600\\
			5	1.521692\\
			6	1.368825\\
			7	1.171699\\
			8	1.010232\\
			9	0.877973\\
			10	0.769640\\
			11	0.680903\\
			12	0.608218\\
			13	0.548682\\
			14	0.499067\\
			15	0.487392\\
			16	0.477939\\
			17	0.471717\\
			18	0.418345\\
			19	0.374570\\
			20	0.338714\\
			21	0.309369\\
			22	0.285363\\
			23	0.265732\\
			24	0.249681\\
			25	0.236557\\
			26	0.225828\\
			27	0.217056\\
			28	0.191916\\
			29	0.170417\\
			30	0.152212\\
			31	0.136918\\
			32	0.124150\\
			33	0.113543\\
			34	0.104766\\
			35	0.097525\\
			36	0.091567\\
			37	0.086675\\
			38	0.082664\\
			39	0.079380\\
			40	0.076695\\
			41	0.074502\\
			42	0.072710\\
			43	0.071249\\
			44	0.070056\\
			45	0.069084\\
			46	0.068292\\
			47	0.067647\\
			48	0.067121\\
			49	0.066693\\
			50	0.066344\\
			51	0.066060\\
			52	0.065829\\
			53	0.065641\\
			54	0.065487\\
			55	0.065363\\
			56	0.065261\\
			57	0.065178\\
			58	0.065111\\
			59	0.065056\\
			60	0.065011\\
			61	0.064975\\
			62	0.064945\\
			63	0.064921\\
			64	0.064902\\
			65	0.064886\\
			66	0.064873\\
			67	0.064863\\
			68	0.064854\\
			69	0.064847\\
			70	0.064842\\
			71	0.064837\\
			72	0.064833\\
			73	0.064830\\
			74	0.064828\\
			75	0.064825\\
			76	0.064824\\
			77	0.064822\\
			78	0.064821\\
			79	0.064820\\
			80	0.064820\\
			81	0.064819\\
			82	0.064818\\
			83	0.064818\\
			84	0.064818\\
			85	0.064817\\
			86	0.064817\\
			87	0.064817\\
			88	0.064817\\
			89	0.064817\\
			90	0.064816\\
			91	0.064816\\
			92	0.064817\\
			93	0.064816\\
			94	0.064817\\
			95	0.064817\\
			96	0.064817\\
			97	0.064817\\
			98	0.064816\\
			99	0.061755\\
			100	0.059251\\
		};
		\addlegendentry{$z_{k,2}$}
	\end{axis}
\end{tikzpicture}
\hspace{-1.3em}\begin{tikzpicture}
\begin{axis}[%
		width=0.3\textwidth,
		xmin=0,
		xmax=100,
		ymin=-5,
		ymax=25,
		grid=both,
		grid style={line width=.1pt, draw=gray!10},
    	minor tick num=1,
		tick label style={font=\footnotesize},
		legend cell align=left,
		legend style={at={(axis cs:98,24.3)},anchor=north east,nodes={scale=0.8, transform shape}}
		]
		\addplot[
color=blue] table[row sep=crcr]{%
1	16.2613872\\
2	12.0703496\\
3	7.56668695\\
4	3.94583570\\
5	2.60449615\\
6	4.95320818\\
7	4.09822459\\
8	1.83262167\\
9	3.29447104\\
10	-1.2347379\\
11	6.20344147\\
12	3.12976431\\
13	1.66415381\\
14	-0.5265442\\
15	2.83635504\\
16	4.26350743\\
17	-1.9805784\\
18	-0.6798347\\
19	0.57826553\\
20	2.78723565\\
21	1.38251118\\
22	2.37101629\\
23	-0.6881477\\
24	0.41686131\\
25	-0.4901132\\
26	2.44690448\\
27	2.50421513\\
28	-0.7823254\\
29	4.20850015\\
30	2.12285294\\
31	0.54585377\\
32	1.31415922\\
33	-0.2799353\\
34	-1.1121442\\
35	-0.5986040\\
36	3.99682967\\
37	-1.0838351\\
38	3.57581346\\
39	1.39912178\\
40	3.52903847\\
41	5.84812033\\
42	3.19768782\\
43	-1.4465213\\
44	2.23016511\\
45	-1.9507211\\
46	-0.3765207\\
47	-1.7793049\\
48	0.00578815\\
49	-3.0622652\\
50	2.78481978\\
51	-0.7767775\\
52	3.36929606\\
53	-2.7862210\\
54	1.23893907\\
55	0.35153424\\
56	0.35462467\\
57	-1.5337390\\
58	-1.5597646\\
59	1.36056189\\
60	3.11787730\\
61	1.67727331\\
62	1.78345550\\
63	0.59435776\\
64	2.86651345\\
65	-0.6132170\\
66	-0.6255361\\
67	1.59298326\\
68	5.99683251\\
69	4.21720170\\
70	4.26144265\\
71	2.92836608\\
72	4.43880920\\
73	-1.1514508\\
74	2.47067618\\
75	-2.2356449\\
76	0.85185322\\
77	-0.0303331\\
78	4.92055466\\
79	3.33570055\\
80	3.30083744\\
81	3.45774644\\
82	2.01535217\\
83	-1.1488658\\
84	-0.7333348\\
85	1.70233158\\
86	3.28832164\\
87	1.85646730\\
88	-2.2393505\\
89	2.72439463\\
90	3.79236316\\
91	-1.1037253\\
92	0.16355904\\
93	3.28276342\\
94	0.88897364\\
95	1.31582572\\
96	2.99603075\\
97	-2.1662616\\
98	2.25744724\\
99	1.61166214\\
100	0.26181649\\
};
\addlegendentry{$\hat{\phi}(z_k)+w_{k,1}$}
\addplot[smooth, color=red] table[row sep=crcr]{%
1	19      \\
2	13.48204\\
3	8.756411\\
4	6.143274\\
5	5.136213\\
6	4.012588\\
7	3.084951\\
8	2.457351\\
9	2.032000\\
10	1.743118\\
11	1.546430\\
12	1.412118\\
13	1.320079\\
14	1.256155\\
15	1.237711\\
16	1.230510\\
17	1.230548\\
18	1.184935\\
19	1.151515\\
20	1.126772\\
21	1.108266\\
22	1.094288\\
23	1.083630\\
24	1.075431\\
25	1.069073\\
26	1.064105\\
27	1.060198\\
28	1.042134\\
29	1.031032\\
30	1.023801\\
31	1.018878\\
32	1.015416\\
33	1.012919\\
34	1.011083\\
35	1.009711\\
36	1.008673\\
37	1.007878\\
38	1.007265\\
39	1.006787\\
40	1.006412\\
41	1.006116\\
42	1.005881\\
43	1.005694\\
44	1.005545\\
45	1.005425\\
46	1.005328\\
47	1.005251\\
48	1.005188\\
49	1.005137\\
50	1.005096\\
51	1.005063\\
52	1.005036\\
53	1.005014\\
54	1.004996\\
55	1.004982\\
56	1.004970\\
57	1.004961\\
58	1.004953\\
59	1.004947\\
60	1.004941\\
61	1.004937\\
62	1.004934\\
63	1.004931\\
64	1.004929\\
65	1.004927\\
66	1.004926\\
67	1.004925\\
68	1.004924\\
69	1.004923\\
70	1.004922\\
71	1.004922\\
72	1.004921\\
73	1.004921\\
74	1.004921\\
75	1.004920\\
76	1.004920\\
77	1.004920\\
78	1.004920\\
79	1.004920\\
80	1.004920\\
81	1.004919\\
82	1.004919\\
83	1.004919\\
84	1.004919\\
85	1.004919\\
86	1.004919\\
87	1.004919\\
88	1.004919\\
89	1.004919\\
90	1.004919\\
91	1.004919\\
92	1.004919\\
93	1.004919\\
94	1.004919\\
95	1.004919\\
96	1.004919\\
97	1.004919\\
98	1.004919\\
99	1.004075\\
100	1.003566\\
};
\addlegendentry{$\hat{\phi}(z_k)$}
\addplot[smooth, color=green!70!black] table[row sep=crcr]{%
1	21.74000000\\
2	16.72915836\\
3	9.476743120\\
4	6.658103725\\
5	6.086521079\\
6	4.727921705\\
7	3.652433373\\
8	3.060962178\\
9	2.664124110\\
10	2.397871994\\
11	2.219234425\\
12	2.099380425\\
13	2.018966318\\
14	1.965783482\\
15	1.838522458\\
16	1.803894757\\
17	1.784378275\\
18	1.497717065\\
19	1.478290283\\
20	1.465370919\\
21	1.456763464\\
22	1.451021729\\
23	1.447188543\\
24	1.444628112\\
25	1.442917231\\
26	1.441773753\\
27	1.441009381\\
28	1.202971247\\
29	1.190483830\\
30	1.183429866\\
31	1.179255024\\
32	1.176697224\\
33	1.175091010\\
34	1.174064797\\
35	1.173401285\\
36	1.172968829\\
37	1.172685446\\
38	1.172499066\\
39	1.172376193\\
40	1.172295055\\
41	1.172241424\\
42	1.172205948\\
43	1.172182470\\
44	1.172166928\\
45	1.172156637\\
46	1.172149822\\
47	1.172145308\\
48	1.172142319\\
49	1.172140339\\
50	1.172139027\\
51	1.172138159\\
52	1.172137584\\
53	1.172137202\\
54	1.172136950\\
55	1.172136782\\
56	1.172136671\\
57	1.172136599\\
58	1.172136550\\
59	1.172136518\\
60	1.172136497\\
61	1.172136482\\
62	1.172136473\\
63	1.172136466\\
64	1.172136462\\
65	1.172136460\\
66	1.172136458\\
67	1.172136457\\
68	1.172136456\\
69	1.172136455\\
70	1.172136455\\
71	1.172136455\\
72	1.172136455\\
73	1.172136455\\
74	1.172136455\\
75	1.172136455\\
76	1.172136455\\
77	1.172136455\\
78	1.172136455\\
79	1.172136454\\
80	1.172136455\\
81	1.172136455\\
82	1.172136455\\
83	1.172136455\\
84	1.172136455\\
85	1.172136454\\
86	1.172136455\\
87	1.172136455\\
88	1.172136455\\
89	1.172136455\\
90	1.172136455\\
91	1.172136455\\
92	1.172136454\\
93	1.172136455\\
94	1.172136455\\
95	1.172136455\\
96	1.172136455\\
97	1.172136455\\
98	1.172136455\\
99	1.166914901\\
100	1.166101834\\
};
\addlegendentry{$\phi^+(z_k;\Gamma_0\!\cap\!\scalebox{0.6}{\ldots}\!\cap\!\Gamma_k)\!\!$}
\end{axis}
\end{tikzpicture}\hspace{-2.5em}%
%
\end{center}
\vspace{-1em}\caption{The simulation results for the initial condition
	$z_0 =\begin{pmatrix} 3 & 3\end{pmatrix}^\top$. In the first plot,
	the resulting trajectory of the estimate $z_k$, shown
	elementwise. The second plot shows the measurements corresponding
	to $z_k+f_1 = z_k$, that is, to $\hat{\phi}(z_k)+w_{k,1}$, where $w_{k,1}$ is
	the first element of $W_k$. These measurements are compared to the
	actual value of the function $\hat{\phi}$ and the current upper bound.} 
  \label{fig:plots single meas}
\end{figure}
%
\begin{figure}[ht]
	\begin{center}
		\vspace{-1em}\hspace{-0.4em}\begin{tikzpicture}
	\begin{axis}[%
		width=0.3\textwidth,
		xmin=0,
		xmax=100,
		ymin=0,
		ymax=3,
		grid=both,
		grid style={line width=.1pt, draw=gray!10},
		minor tick num=1,
		tick label style={font=\footnotesize},
		legend style={at={(axis cs:95,2.8)},anchor=north east,nodes={scale=0.8, transform shape}}
		]
		\addplot[smooth, color=blue] table[row sep=crcr]{%
1	3      \\
2	2.49820\\
3	2.21066\\
4	1.95975\\
5	1.73516\\
6	1.53384\\
7	1.35411\\
8	1.19459\\
9	1.05386\\
10	0.93037\\
11	0.82251\\
12	0.72867\\
13	0.64731\\
14	0.57694\\
15	0.51622\\
16	0.46393\\
17	0.41897\\
18	0.38036\\
19	0.34724\\
20	0.31885\\
21	0.29454\\
22	0.27374\\
23	0.25595\\
24	0.24073\\
25	0.22773\\
26	0.21662\\
27	0.20714\\
28	0.19904\\
29	0.19212\\
30	0.18622\\
31	0.18118\\
32	0.17688\\
33	0.17321\\
34	0.17009\\
35	0.16742\\
36	0.16514\\
37	0.16320\\
38	0.16154\\
39	0.16013\\
40	0.15892\\
41	0.15789\\
42	0.15702\\
43	0.15627\\
44	0.15563\\
45	0.15509\\
46	0.15463\\
47	0.15423\\
48	0.15389\\
49	0.15361\\
50	0.15336\\
51	0.15315\\
52	0.15297\\
53	0.15282\\
54	0.15269\\
55	0.15258\\
56	0.15249\\
57	0.15241\\
58	0.15234\\
59	0.15228\\
60	0.15223\\
61	0.15219\\
62	0.15215\\
63	0.15212\\
64	0.15209\\
65	0.15207\\
66	0.15205\\
67	0.15204\\
68	0.15202\\
69	0.15201\\
70	0.15200\\
71	0.15199\\
72	0.15198\\
73	0.15198\\
74	0.15197\\
75	0.15197\\
76	0.15196\\
77	0.15196\\
78	0.15196\\
79	0.15196\\
80	0.15195\\
81	0.15195\\
82	0.15195\\
83	0.15195\\
84	0.15195\\
85	0.15195\\
86	0.15195\\
87	0.15195\\
88	0.15195\\
89	0.15195\\
90	0.15195\\
91	0.15194\\
92	0.15194\\
93	0.15194\\
94	0.15194\\
95	0.15194\\
96	0.15194\\
97	0.15194\\
98	0.15194\\
99	0.15194\\
100	0.15194\\
		};
		\addlegendentry{$z_{k,1}$}
		\addplot[smooth, color=red] table[row sep=crcr]{%
1	3     \\
2	2.4982\\
3	2.2106\\
4	1.9597\\
5	1.7351\\
6	1.5338\\
7	1.3541\\
8	1.1945\\
9	1.0538\\
10	0.9303\\
11	0.8225\\
12	0.7286\\
13	0.6473\\
14	0.5769\\
15	0.5162\\
16	0.4639\\
17	0.4189\\
18	0.3803\\
19	0.3472\\
20	0.3188\\
21	0.2945\\
22	0.2737\\
23	0.2559\\
24	0.2407\\
25	0.2277\\
26	0.2166\\
27	0.2071\\
28	0.1990\\
29	0.1921\\
30	0.1862\\
31	0.1811\\
32	0.1768\\
33	0.1732\\
34	0.1700\\
35	0.1674\\
36	0.1651\\
37	0.1632\\
38	0.1615\\
39	0.1601\\
40	0.1589\\
41	0.1578\\
42	0.1570\\
43	0.1562\\
44	0.1556\\
45	0.1550\\
46	0.1546\\
47	0.1542\\
48	0.1538\\
49	0.1536\\
50	0.1533\\
51	0.1531\\
52	0.1529\\
53	0.1528\\
54	0.1526\\
55	0.1525\\
56	0.1524\\
57	0.1524\\
58	0.1523\\
59	0.1522\\
60	0.1522\\
61	0.1521\\
62	0.1521\\
63	0.1521\\
64	0.1520\\
65	0.1520\\
66	0.1520\\
67	0.1520\\
68	0.1520\\
69	0.1520\\
70	0.1520\\
71	0.1519\\
72	0.1519\\
73	0.1519\\
74	0.1519\\
75	0.1519\\
76	0.1519\\
77	0.1519\\
78	0.1519\\
79	0.1519\\
80	0.1519\\
81	0.1519\\
82	0.1519\\
83	0.1519\\
84	0.1519\\
85	0.1519\\
86	0.1519\\
87	0.1519\\
88	0.1519\\
89	0.1519\\
90	0.1519\\
91	0.1519\\
92	0.1519\\
93	0.1519\\
94	0.1519\\
95	0.1519\\
96	0.1519\\
97	0.1519\\
98	0.1519\\
99	0.1519\\
100	0.1519\\
		};
		\addlegendentry{$z_{k,2}$}
	\end{axis}
\end{tikzpicture}
\hspace{-1.3em}\begin{tikzpicture}
\begin{axis}[%
		width=0.3\textwidth,
		xmin=0,
		xmax=100,
		ymin=0,
		ymax=22.5,
		grid=both,
		grid style={line width=.1pt, draw=gray!10},
		minor tick num=1,
		ytick ={0,5,10,15,20},
		tick label style={font=\footnotesize},
		legend cell align=left,
		legend style={at={(axis cs:98,22)},anchor=north east,nodes={scale=0.8, transform shape}}
		]
		\addplot[
color=blue] table[row sep=crcr]{%
1	16.2613\\
2	16.2206\\
3	13.5127\\
4	11.4198\\
5	9.76021\\
6	8.44398\\
7	7.40586\\
8	6.59273\\
9	5.95987\\
10	5.46980\\
11	5.09167\\
12	4.80055\\
13	4.57663\\
14	4.40433\\
15	4.27159\\
16	4.16909\\
17	4.08969\\
18	4.02796\\
19	3.97976\\
20	3.94195\\
21	3.91212\\
22	3.88848\\
23	3.86963\\
24	3.85452\\
25	3.84234\\
26	3.83246\\
27	3.82442\\
28	3.81784\\
29	3.81243\\
30	3.80797\\
31	3.80426\\
32	3.80119\\
33	3.79862\\
34	3.79647\\
35	3.79467\\
36	3.79315\\
37	3.79188\\
38	3.79080\\
39	3.78989\\
40	3.78912\\
41	3.78847\\
42	3.78792\\
43	3.78745\\
44	3.78705\\
45	3.78672\\
46	3.78643\\
47	3.78618\\
48	3.78598\\
49	3.78580\\
50	3.78565\\
51	3.78552\\
52	3.78541\\
53	3.78532\\
54	3.78524\\
55	3.78517\\
56	3.78512\\
57	3.78507\\
58	3.78502\\
59	3.78499\\
60	3.78496\\
61	3.78493\\
62	3.78491\\
63	3.78489\\
64	3.78488\\
65	3.78486\\
66	3.78485\\
67	3.78484\\
68	3.78483\\
69	3.78482\\
70	3.78482\\
71	3.78481\\
72	3.78481\\
73	3.78481\\
74	3.78480\\
75	3.78480\\
76	3.78480\\
77	3.78480\\
78	3.78479\\
79	3.78479\\
80	3.78479\\
81	3.78479\\
82	3.78479\\
83	3.78479\\
84	3.78479\\
85	3.78479\\
86	3.78479\\
87	3.78479\\
88	3.78479\\
89	3.78479\\
90	3.78479\\
91	3.78479\\
92	3.78479\\
93	3.78479\\
94	3.78478\\
95	3.78478\\
96	3.78478\\
97	3.78478\\
98	3.78478\\
99	3.78478\\
100	3.78478\\
};
\addlegendentry{$\hat{\phi}(z_k)+w_{k,1}$}
\addplot[smooth, color=red] table[row sep=crcr]{%
1	19      \\
2	13.48204\\
3	10.77409\\
4	8.681271\\
5	7.021605\\
6	5.705369\\
7	4.667254\\
8	3.854124\\
9	3.221259\\
10	2.731187\\
11	2.353059\\
12	2.061945\\
13	1.838020\\
14	1.665726\\
15	1.532982\\
16	1.430477\\
17	1.351081\\
18	1.289353\\
19	1.241155\\
20	1.203338\\
21	1.173516\\
22	1.149871\\
23	1.131021\\
24	1.115909\\
25	1.103727\\
26	1.093856\\
27	1.085816\\
28	1.079235\\
29	1.073825\\
30	1.069358\\
31	1.065656\\
32	1.062577\\
33	1.060009\\
34	1.057861\\
35	1.056059\\
36	1.054544\\
37	1.053269\\
38	1.052193\\
39	1.051284\\
40	1.050515\\
41	1.049864\\
42	1.049312\\
43	1.048843\\
44	1.048446\\
45	1.048108\\
46	1.047821\\
47	1.047577\\
48	1.047369\\
49	1.047192\\
50	1.047041\\
51	1.046913\\
52	1.046804\\
53	1.046711\\
54	1.046632\\
55	1.046564\\
56	1.046507\\
57	1.046458\\
58	1.046416\\
59	1.046380\\
60	1.046350\\
61	1.046324\\
62	1.046302\\
63	1.046283\\
64	1.046267\\
65	1.046253\\
66	1.046242\\
67	1.046232\\
68	1.046224\\
69	1.046217\\
70	1.046211\\
71	1.046205\\
72	1.046201\\
73	1.046197\\
74	1.046194\\
75	1.046191\\
76	1.046189\\
77	1.046187\\
78	1.046185\\
79	1.046184\\
80	1.046183\\
81	1.046182\\
82	1.046181\\
83	1.046180\\
84	1.046180\\
85	1.046179\\
86	1.046179\\
87	1.046178\\
88	1.046178\\
89	1.046177\\
90	1.046177\\
91	1.046177\\
92	1.046177\\
93	1.046177\\
94	1.046177\\
95	1.046177\\
96	1.046177\\
97	1.046176\\
98	1.046176\\
99	1.046176\\
100	1.046176\\
};
\addlegendentry{$\hat{\phi}(z_k)$}
\addplot[smooth, color=green!70!black] table[row sep=crcr]{%
1	21.7400000\\
2	16.7291583\\
3	15.0963846\\
4	13.5824812\\
5	12.2340045\\
6	11.0883293\\
7	10.1435102\\
8	9.38106812\\
9	8.77549187\\
10	8.29991876\\
11	7.92933520\\
12	7.64201180\\
13	7.41987708\\
14	7.24833038\\
15	7.11580811\\
16	7.01327783\\
17	6.93375300\\
18	6.87186677\\
19	6.82351460\\
20	6.78556402\\
21	6.75562968\\
22	6.73189541\\
23	6.71297525\\
24	6.69780990\\
25	6.68558799\\
26	6.67568603\\
27	6.66762315\\
28	6.66102593\\
29	6.65560358\\
30	6.65112794\\
31	6.64741920\\
32	6.64433594\\
33	6.64176427\\
34	6.63961343\\
35	6.63780985\\
36	6.63629418\\
37	6.63501772\\
38	6.63394120\\
39	6.63303180\\
40	6.63226259\\
41	6.63161112\\
42	6.63105884\\
43	6.63059052\\
44	6.63019289\\
45	6.62985511\\
46	6.62956801\\
47	6.62932380\\
48	6.62911595\\
49	6.62893902\\
50	6.62878843\\
51	6.62866020\\
52	6.62855096\\
53	6.62845792\\
54	6.62837880\\
55	6.62831147\\
56	6.62825414\\
57	6.62820529\\
58	6.62816350\\
59	6.62812774\\
60	6.62809725\\
61	6.62807143\\
62	6.62804934\\
63	6.62803069\\
64	6.62801450\\
65	6.62800075\\
66	6.62798902\\
67	6.62797928\\
68	6.62797096\\
69	6.62796385\\
70	6.62795779\\
71	6.62795262\\
72	6.62794821\\
73	6.62794444\\
74	6.62794123\\
75	6.62793849\\
76	6.62793616\\
77	6.62793420\\
78	6.62793252\\
79	6.62793110\\
80	6.62792988\\
81	6.62792885\\
82	6.62792798\\
83	6.62792725\\
84	6.62792665\\
85	6.62792614\\
86	6.62792565\\
87	6.62792526\\
88	6.62792480\\
89	6.62792455\\
90	6.62792436\\
91	6.62792419\\
92	6.62792407\\
93	6.62792395\\
94	6.62792384\\
95	6.62792370\\
96	6.62792362\\
97	6.62792357\\
98	6.62792355\\
99	6.62792354\\
100	6.62792354\\
};
\addlegendentry{$\phi^+(z_k;\Gamma_0\!\cap\!\scalebox{0.6}{\ldots}\!\cap\!\Gamma_k)\!\!$}
\end{axis}
\end{tikzpicture}\hspace{-2.5em}%
%
	\end{center}
\vspace{-1em}  \caption{The simulation results for the initial condition
	$z_0 =\begin{pmatrix} 3 & 3\end{pmatrix}^\top$ with nonrandom
	noise. In particular, the noise samples are equal and proportional
	to the vector of ones. The figures correspond to those in
	Figure~\ref{fig:plots single meas}. Note that the elements of the
	vector $z$ are equal for all time. }
	\vspace{-1.5em}
\label{fig:nonrandom}
\end{figure}
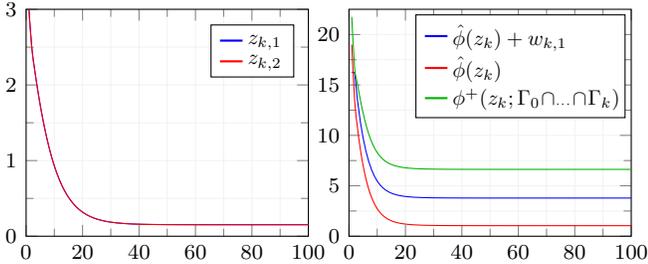

Lastly, we illustrate that worst-case, adversarial noise can lead to
convergence to a suboptimal bound. For this, in the last set of
simulations and for the same scenario, instead of generating noise
randomly, we apply the same noise sample from $\calZ(\Pi)$ at each
step after $k=1$. These results can be seen in
Figure~\ref{fig:nonrandom}.  In particular, note that $z_k$ does not
converge to the optimizer at the origin origin, but to the point
$\begin{pmatrix} 0.1519 &0.1519\end{pmatrix}^\top$. Moreover, since
the noise is constant, it can be seen that the upper bound does not
converge to the optimal value of the true function.

\section{Conclusions}
We have investigated suboptimization for unknown functions on the
basis of measurements with bounded noise. Employing ideas from the
informativity framework for data-driven control, the notions of
set-valued regression and the cautious suboptimization problem were
introduced. In short, the data gives rise to a set of possible 
parameters, and in order to draw conclusions regarding the true 
function we require bounds for all possible realizations of this parameter. 
Resolving this problem was shown to be equivalent to minimization of the
worst-case realization. For this, we provided explicit forms and 
convexity results, allowing efficient solutions. In an online
setting, we investigated the iteration of cautious suboptimization and
local collection of new measurements. This procedure gives rise to 
nonincreasing guaranteed upper bounds for the optimal value of the unknown
function. Moreover, in the case that the noise is randomly generated,
this procedure is proven to converge to the true optimal value.

A number of avenues for future work present themselves.  As
illustrated by set membership estimation (see e.g. \cite{MM-AV:91}),
the Lipschitz constant is a powerful tool for deriving local bounds on
the basis of measurements. Indeed, Lipschitz constants of the
parameterized functions can be derived from those of the basis
functions, allowing for the determination of locally suboptimal values
outside the scope of convex functions. Another extension would be to
consider non-scalar functions, with an aim at analysis of nonlinear
systems within a parameterized class. We also would like to
characterize the sample efficiency and computational complexity of the
proposed techniques.  Regarding online methods, there are many
possible extensions, including relaxations of the optimization problem
to avoid intersections of a large amount of convex sets. This work
investigates uniformly distributed noise, and investigation of
different distributions and less conservative convergence results is a
topic of interest.

\bibliography{../bib/alias,../bib/JC,../bib/Main,../bib/Main-add,../bib/New,../bib/FB}
\bibliographystyle{IEEEtran}

\end{document}